# Weak convergence of error processes in discretizations of stochastic integrals and Besov spaces

STEFAN GEISS* and ANNI TOIVOLA**

*Department of Mathematics and Statistics, P.O. Box 35, FIN-40014 University of Jyväskylä, Finland. E-mails: \*geiss@maths.jyu.fi; \*\*atoivola@maths.jyu.fi*

We consider weak convergence of the rescaled error processes arising from Riemann discretizations of certain stochastic integrals and relate the $L_p$-integrability of the weak limit to the fractional smoothness in the Malliavin sense of the stochastic integral.

## 1. Introduction

Quantitative approximation problems for stochastic integrals appear naturally in stochastic finance. Consider a stochastic integral

$$g(X_1) = \mathbb{E}g(X_1) + \int_0^1 \frac{\partial G}{\partial x}(u, X_u) \, dX_u,$$

where the diffusion $X = (X_t)_{t \in [0,1]}$ models a price process, $g(X_1) \in L_2$ is a pay-off of a European type option, and $G$ solves a corresponding parabolic backward partial differential equation (PDE) with $g(x) = G(1,x)$. We look at the Riemann approximation

$$\sum_{i=0}^{n-1} \frac{\partial G}{\partial x}(t_i, X_{t_i})(X_{t_{i+1}} - X_{t_i})$$

along a deterministic time-net $\tau = (t_i)_{i=1}^n$ and the error process $C(\tau) = (C_t(\tau))_{t \in [0,1]}$ given by

$$C_t(\tau) := \int_0^t \frac{\partial G}{\partial x}(u, X_u) \, dX_u - \sum_{i=0}^{n-1} \frac{\partial G}{\partial x}(t_i, X_{t_i})(X_{t_{i+1} \wedge t} - X_{t_i \wedge t}).$$

The process $C(\tau)$ describes the hedging error that occurs when a continuously adjusted portfolio is replaced by a portfolio that is adjusted only at the time-knots $t_0, \ldots, t_{n-1}$. Given a sequence of time-nets $\tau^n = (t_i^n)_{i=0}^n$, one is interested in the rate of convergence of $C(\tau^n)$ towards zero as $n \to \infty$. There are (at least) two principal ways to measure the size of $C(\tau^n)$. First, one can use strong criteria, like $L_p$-norms, where one typically looks for estimates of the form

$$\|C_1(\tau^n)\|_{L_p} \leq cn^{-\theta} \tag{1}$$

for some $\theta > 0$. Second, one can investigate the weak convergence of the re-scaled processes $\sqrt{n}C(\tau^n)$. A priori, in our setting there are no general principles known so far to deduce a certain strong convergence from a weak limit or to go the other way around. Concepts of weak convergence are of particular interest in applications because they already provide the needed information in many cases and promise potentially better approximation rates than obtained under strong criteria. Results about weak convergence in our context are obtained for example in [22] (see also [4]) and [13, 15, 24]. For the general theory the reader is referred to [18] (see also [21]).

Gobet and Temam have shown in [13], Theorems 1 and 3, that for the binary option (i.e., $g(x) = \chi_{[K,\infty)}(x)$ for $K > 0$ and $X$ being the geometric Brownian motion), in case of equidistant time-nets $\tau^n$, the scaling factor for the weak convergence can be taken to be $n^{1/2}$ whereas the $L_2$-rate in (1) is $\theta = \frac{1}{4}$. Intuitively, one would expect that the scaling exponent $\frac{1}{2}$ and $\theta$ coincide. Indeed, for pay-off functions $g$ having some fractional smoothness in the Malliavin sense (like the binary option with smoothness $\beta \in (0, \frac{1}{2})$; see Example 2.2), the $L_2$-rate $\theta = \frac{1}{2}$ can always be achieved by using appropriate non-equidistant time-nets, see [6, 12]. Non-equidistant time-nets have been also used in other papers like, for example, [14, 20].

From this, two questions become natural: Is there a connection between fractional smoothness and weak convergence? And, do non-equidistant time-nets have a positive effect on the weak convergence?

The aim of this paper is to answer both questions in the positive at the same time by investigating the $L_p$-integrability of the weak limit of $\sqrt{n}C(\tau^n)$ as $n \to \infty$ for different sequences of time-nets $\tau^n$. This has relevance for applications where good tail estimates for the weak limit are desirable.

The paper is organized as follows:

- After introducing the notation we formulate our basic result, Theorem 3.1, where we characterize the existence of a square integrable weak limit of $\sqrt{n}C(\tau^{n,\beta})$ by the condition that $g$ or $g(\exp(\cdot - \frac{1}{2}))$ (depending on the diffusion $X$) belongs to the Besov space $B_{2,2}^\beta(\gamma)$. The parameter $\beta \in (0,1]$ is the fractional smoothness in the Malliavin sense and the time-nets

$$\tau^{n,\beta} := (1 - (1 - (i/n))^{1/\beta})_{i=0}^n$$

are adapted to the smoothness $\beta$. Hence, if $g$ or $g(\exp(\cdot - \frac{1}{2}))$ have a non-trivial fractional smoothness and if we use the right time-nets, then we always get a square-integrable weak limit. The concept of fractional smoothness allows us to consider at once the large class of functions $\bigcup_{\beta \in (0,1]} B_{2,2}^\beta(\gamma)$, which contains all examples usually studied in the literature in this context. For the binary option this means that the weak limit for equidistant time-nets in [13] is not square-integrable but becomes square-integrable for the time-nets $\tau^{n,\beta}$ as long as $\beta \in (0, \frac{1}{2})$ because of Example 2.2 below.

- The $L_2$-setting of Theorem 3.1 is extended in Section 4 to the $L_p$-setting, $p \in [2, \infty)$. Corollary 4.4 gives nearly optimal conditions that the weak limit is $L_p$-integrable. As an application for the binary option we compute in Example 4.7 the best possible $L_p$-integrability of the weak limit provided that the $\tau^{n,\beta}$-nets are used. In particular, this example shows how the integrability can be improved to any $p \in [2, \infty)$ by using nets $\tau^{n,\beta}$ with a $\beta$ small enough or, equivalently, by using nets with a sufficiently high concentration of the time-knots close to the final time-point $t = 1$.
- The upper estimate for the $L_p$-integrability of the weak limit in Example 4.7 for the binary option has a more general background: in Corollary 4.11 we assume that $g$ has a local singularity of order $\eta \geq 0$ (measured in terms of a sharp function) and deduce an upper bound for the $L_p$-integrability of the weak limit.
- In Section 5 we prove Theorem 3.1. First, we derive the existence of the weak limit on $[0, T]$ with $T \in (0, 1)$. Second, as the main part, we have to deal with a singularity of our approximation problem at time $t = 1$ because of the blow-up of the Malliavin derivative of $\mathbb{E}(g(X_1)|\mathcal{F}_t)$ as $t \uparrow 1$. The degree of this blow-up is connected to the fractional smoothness $\beta$ of $g$. The used time-nets $\tau^{n,\beta}$ are essential as they are chosen to compensate this singularity.

In this paper we restrict ourselves to the case that the driving process of the stochastic integrals is the Brownian motion or the geometric Brownian motion. This is done to connect weak convergence and Besov spaces as exactly as possible where we use results from [12] proved for this setting. A setting for Lévy–Itô processes was considered in [23] and a setting under transaction costs in [5]. Extensions along the line of fractional smoothness and non-homogeneous time-nets might be investigated in the following directions:

- To consider more general diffusions as driving process $X$, the setting of [6] seems to be appropriate for a first step where the scale of Besov spaces was replaced by a scale that has more appropriate stochastic descriptions and where one imposes quantitative smoothness assumptions on the parameters of the diffusion. This extension would still result in a setting where an underlying PDE exists so that the proofs are expected to be parallel to the ones of this paper. Note that our integrand $(\partial G / \partial x)(u, X_u)$ is obtained by the PDE (3) below.
- One might also investigate stochastic integrals where the integrands are not obtained via a PDE. In this case appropriate structural assumptions should be necessary as in the present paper the blow-up of the integrands at time 1 is essential and the right notion of fractional smoothness has to be taken.

There are other settings where weak limits are investigated for rescaled processes that have a similar structure to $\sqrt{n} C(\tau^{n,\beta})$. Such a situation is the consideration of weak limits of rescaled error processes arising in Euler schemes for diffusions, see [17]. It would be of interest to check further extensions of the present paper using these results and ideas.

## 2. Notation

Let $B = (B_t)_{t\in[0,1]}$ be a standard Brownian motion defined on a stochastic basis $(\Omega, \mathcal{F}, \mathbb{P}, (\mathcal{F}_t)_{t\in[0,1]})$, where $B_0 \equiv 0$, all paths are assumed to be continuous, $(\mathcal{F}_t)_{t\in[0,1]}$ is the augmentation of the natural filtration of $B$ and $\mathcal{F} = \mathcal{F}_1$. Let $X$ be either the Brownian motion or the geometric Brownian motion $S = (S_t)_{t\in[0,1]}$ with

$$S_t := e^{B_t - t/2}.$$

To treat both cases for $X$ simultaneously we let $\sigma(x) \equiv 1$ if $X = B$ and $\sigma(x) = x$ if $X = S$, so that $dX_t = \sigma(X_t)\, dB_t$. Let $g : E \to \mathbb{R}$ be a Borel function, where $E = \mathbb{R}$ if $X = B$ and $E = (0, \infty)$ if $X = S$, such that $\mathbb{E}g(X_1)^2 < \infty$. Define the function $G$ by setting

$$G(t, x) := \mathbb{E}(g(X_1) | X_t = x) = \begin{cases} \mathbb{E}g(x + X_{1-t}), & X = B, \\ \mathbb{E}g(x X_{1-t}), & X = S. \end{cases} \quad (2)$$

Then it follows that $G \in C^\infty([0,1) \times E)$ and that $G$ satisfies the partial differential equation

$$\frac{\partial G}{\partial t}(t, x) + \frac{\sigma(x)^2}{2} \frac{\partial^2 G}{\partial x^2}(t, x) = 0 \quad (3)$$

for $(t, x) \in [0, 1) \times E$ with $G(1, x) = g(x)$. By Itô's formula,

$$g(X_1) = \mathbb{E}g(X_1) + \int_0^1 \frac{\partial G}{\partial x}(u, X_u)\, dX_u \quad \text{a.s.}$$

Our interest is to approximate the stochastic integral $\int_0^1 \frac{\partial G}{\partial x}(u, X_u)\, dX_u$ by a Riemann approximation. To this end, given a deterministic time-net $\tau = (t_i)_{i=0}^n$ with $0 = t_0 < \cdots < t_n = 1$, we define the error process $C(\tau) = (C_t(\tau))_{t\in[0,1]}$ by

$$C_t(\tau) := \int_0^t \frac{\partial G}{\partial x}(u, X_u)\, dX_u - \sum_{i=0}^{n-1} \frac{\partial G}{\partial x}(t_i, X_{t_i})(X_{t_{i+1}\wedge t} - X_{t_i \wedge t})$$

for $t \in [0, 1]$, where we can assume that all paths are continuous. For $\beta \in (0, 1]$ we introduce the time-nets $\tau^{n,\beta} := (t_i^{n,\beta})_{i=0}^n$ defined by

$$t_i^{n,\beta} := 1 - \left(1 - \frac{i}{n}\right)^{1/\beta}.$$

The smaller the $\beta$, the higher the concentration of the time-knots is near to one. In particular,

$$\frac{|t^{n,\beta}_{i+1} - u|}{(1-u)^{1-\beta}} \le \frac{|t^{n,\beta}_{i+1} - t^{n,\beta}_i|}{(1-t^{n,\beta}_i)^{1-\beta}} \le \frac{1}{\beta n} \qquad \text{for } u \in [t^{n,\beta}_i, t^{n,\beta}_{i+1}) \tag{4}$$

and all $n = 1, 2, \ldots$ and $i = 0, \ldots, n-1$. The Besov spaces we use can be described by Hermite expansions as follows:

**Definition 2.1.** *Let $\mathrm{d}\gamma(x) = (1/\sqrt{2\pi}) \exp(-x^2/2) \, \mathrm{d}x$ be the standard Gaussian measure on $\mathbb{R}$ and let $(h_k)_{k=0}^\infty \subset L_2(\gamma)$ be the orthonormal basis consisting of Hermite polynomials obtained by*

$$h_k(x) := \frac{(-1)^k}{\sqrt{k!}} \mathrm{e}^{x^2/2} \frac{\mathrm{d}^k}{\mathrm{d}x^k}(\mathrm{e}^{-x^2/2}).$$

*Given $\beta \in (0,1]$ and $f = \sum_{k=0}^\infty \alpha_k h_k$, we let $f \in B^\beta_{2,2}(\gamma)$ provided that*

$$\|f\|_{B^\beta_{2,2}(\gamma)} := \left( \sum_{k=0}^\infty (k+1)^\beta \alpha_k^2 \right)^{1/2} < \infty.$$

The parameter $\beta$ is the degree of fractional smoothness. In particular, we have that $B^1_{2,2}(\gamma)$ is the Malliavin Sobolev space $D_{1,2}(\gamma)$. For $\beta \in (0,1)$ the above Besov spaces can also be obtained by the real interpolation method as

$$B^\beta_{2,2}(\gamma) = (L_2(\gamma), D_{1,2}(\gamma))_{\beta,2},$$

see [2], Theorem 5.6.1. On the other hand, for $\beta \in (0,1]$, $f := g$ if $X = B$ and $f(x) := g(\mathrm{e}^{x-1/2})$ if $X = S$, one has that

$$\|g(X_1)\|^2_{L_2} + \int_0^1 (1-u)^{1-\beta} \mathbb{E}\left| \left( \sigma^2 \frac{\partial^2 G}{\partial x^2} \right)(u, X_u) \right|^2 \mathrm{d}u < \infty \tag{5}$$

if and only if $f \in B^\beta_{2,2}(\gamma)$ and

$$\|g(X_1)\|^2_{L_2} + \sup_{t \in [0,1)} (1-t)^{1-\beta} \int_0^t \mathbb{E}\left| \left( \sigma^2 \frac{\partial^2 G}{\partial x^2} \right)(u, X_u) \right|^2 \mathrm{d}u < \infty \tag{6}$$

if and only if $f \in (L_2(\gamma), D_{1,2}(\gamma))_{\beta,\infty}$. The equivalence (5) follows from [12], proof of Theorem 3.2, and the equivalence (6) is implicitly contained in [12] as well (see the preprint version of [11], proof of Theorem 1.3). Because from the general interpolation theory it is known that

$$(L_2(\gamma), D_{1,2}(\gamma))_{\beta',\infty} \subseteq B^\beta_{2,2}(\gamma) \subseteq (L_2(\gamma), D_{1,2}(\gamma))_{\beta,\infty} \tag{7}$$

for $0 < \beta < \beta' < 1$ the reader can deduce from relations (5)–(7), and from [6], page 358, explicit examples of $f \in B_{2,2}^\beta(\gamma)$. For the binary option mentioned in the introduction we have:

**Example 2.2.** For $g(x) = \chi_{[K,\infty)}(x)$ with $K \in \mathbb{R}$ we have that

$$g \in (L_2(\gamma), D_{1,2}(\gamma))_{1/2,\infty} \subseteq \bigcap_{\beta \in (0,1/2)} B_{2,2}^\beta(\gamma).$$

To formulate our results, given $\beta \in (0,1]$ and $t \in [0,1)$, we define

$$\nu_\beta(t) := \frac{1}{\beta}(1-t)^{1-\beta},$$

$$A_\beta(t) := \frac{1}{2}\int_0^t \nu_\beta(u)\left[\left(\sigma^2 \frac{\partial^2 G}{\partial x^2}\right)(u, X_u)\right]^2 du,$$

$$Z_\beta(t) := W_{A_\beta(t)},$$

where $G$ is obtained from the function $g$ as in (2) and $W = (W_t)_{t \geq 0}$ is a standard Brownian motion starting at zero defined on some auxiliary probability space $(M, \mu)$, where we may and do assume that all paths are continuous. Finally, we extend the process $A_\beta$ by

$$A_\beta(1) := \lim_{t \uparrow 1} A_\beta(t),$$

which might be an extended random variable. In the following $\Longrightarrow_{C[0,T]}$ stands for the weak convergence in $C[0,T]$ for some $T > 0$.

## 3. The basic result

The basic result of the paper is:

**Theorem 3.1.** *Let $\beta \in (0,1]$ and $g(X_1) \in L_2$. Then, for all $T \in [0,1)$,*

$$(\sqrt{n} C_t(\tau^{n,\beta}))_{t \in [0,T]} \Longrightarrow_{C[0,T]} (Z_\beta(t))_{t \in [0,T]} \quad \text{as } n \to \infty. \tag{8}$$

*Moreover, the following assertions are equivalent:*

(i) *One has $g \in B_{2,2}^\beta(\gamma)$ for $X = B$ and $g(e^{\cdot -1/2}) \in B_{2,2}^\beta(\gamma)$ for $X = S$, respectively.*
(ii) *On some stochastic basis there exists a continuous square-integrable martingale $M = (M_t)_{t \in [0,1]}$ such that $\sqrt{n} C(\tau^{n,\beta}) \Longrightarrow_{C[0,1]} M$.*
(iii) *One has $\mathbb{E} A_\beta(1) < \infty$ and for $\widetilde{Z}_\beta(t) := W_{A_\beta(t)\chi_{\{A_\beta(1)<\infty\}}}$ with $t \in [0,1]$ it holds that*

$$\sqrt{n} C(\tau^{n,\beta}) \Longrightarrow_{C[0,1]} (\widetilde{Z}_\beta(t))_{t \in [0,1]}.$$

It should be noted that we do not assume any a priori smoothness assumptions for $g$, only the integrability $g(X_1) \in L_2$. Moreover, for $g \in B_{2,2}^{\beta}(\gamma) \setminus B_{2,2}^{\beta'}(\gamma)$ with $0 < \beta < \beta' \leq 1$ one has that

$$\sup_{t \in [0,1)} \mathbb{E}|Z_{\beta'}(t)|^2 = \infty,$$

which follows directly from Section 5.1. Hence the $L_2$-boundedness of $M$ and $(\widetilde{Z}_\beta(t))_{t \in [0,1]}$ in Theorem 3.1 is due to the proper choice of the time-nets.

As a consequence of Theorem 3.1 we obtain that the weak convergence and the $L_2$-boundedness of the rescaled error processes $\sqrt{n}C(\tau^{n,\beta})$ imply each other:

**Corollary 3.2.** *For $\beta \in (0,1]$ and $g(X_1) \in L_2$ the following assertions are equivalent:*

(i) *One has $\sup_{n \geq 1} \sqrt{n}\|C_1(\tau^{n,\beta})\|_{L_2} < \infty$.*
(ii) *On some stochastic basis there exists a continuous square-integrable martingale $M = (M_t)_{t \in [0,1]}$ such that $\sqrt{n}C(\tau^{n,\beta}) \Longrightarrow_{C[0,1]} M$.*

**Proof.** In [12], Theorem 3.2, it was shown that Theorem 3.1(i) is equivalent to $\sup_{n \geq 1} \sqrt{n}\|C_1(\tau^{n,\beta})\|_{L_2} < \infty$ so that we are done. □

As usual, having a weak convergence like in Theorem 3.1(iii), one obtains the weak convergence of functionals $\varphi(\sqrt{n}C(\tau^{n,\beta}))$ whenever $\varphi: C[0,1] \to \mathbb{R}$ is continuous.

## 4. $L_p$-integrability of the weak limit

Let $\beta \in (0,1]$, $X = S$, and $g(e^{\cdot -1/2}) \in B_{2,2}^\beta(\gamma)$. Assume that

$$(\mathbb{P} \times \mu)(|W_{A_\beta(1)\chi_{\{A_\beta(1) < \infty\}}}| > \lambda) \leq \psi(\lambda)$$

for $\lambda > 0$, where $(M, \mu)$ is the auxiliary probability space on which the independent Brownian motion $W$ is defined and $\psi: [\lambda_0, \infty) \to (0,1]$ is a decreasing bijection for some $\lambda_0 \in [0, \infty)$ extended to $[0, \lambda_0]$ by $\psi(\lambda) \equiv 1$. Then

$$\lim_n \mathbb{P}(\psi^{-1}(\varepsilon)n^{-1/2} + C_1(\tau^{n,\beta}) \leq 0) \leq \varepsilon \qquad \text{for } \varepsilon \in (0,1),$$

where Theorem 3.1 and the fact $(\mathbb{P} \times \mu)(|W_{A_\beta(1)\chi_{\{A_\beta(1) < \infty\}}}| = \lambda) = 0$ for $\lambda > 0$ guarantee that the limit exists.

This means that, considering a European option with pay-off $g$ in the discounted Black–Scholes model, an increase of the initial capital by $\psi^{-1}(\varepsilon)n^{-1/2}$ gives, for large $n$, approximately a shortfall probability of at most $\varepsilon$ if the portfolio is re-balanced along the time-net $\tau^{n,\beta}$. Therefore, in order to minimize the increase of the initial capital to reach the pre-given shortfall probability $\varepsilon$, we have to find functions $\psi$ that decrease as fast as possible. Starting with this motivation we proceed as follows in this section, where we consider both cases $X = B$ and $X = S$ if not stated otherwise.

In Corollary 4.4(i) and Proposition 4.5(iii) we give verifiable conditions that the weak limit in Theorem 3.1 at time $T = 1$ has a $p$th moment, $p \in [2, \infty)$. In our context verifiable means that Corollary 4.4(i) can be checked by solving the PDE (3) (see [6]) and that Proposition 4.5(iii) follows for Hölder continuous functions or from general upper bounds for functions $h$ of bounded variation (see [1], Theorem 2.4). This is connected to fractional smoothness in terms of Besov spaces, which follows from the equivalences in Proposition 4.5 and Remark 4.6. As an application we demonstrate for the pay-off of the binary option in Example 4.7 that the density of the time-knots of our time-nets $\tau^{n,\beta}$ close to maturity directly affects the integrability of the weak limit. Corollary 4.11 accompanies this statement by showing how local properties of $g$ change for the worse the behavior of the error processes $C(\tau^{n,\beta})$ by proving an upper bound for the $L_p$-integrability of their rescaled weak limit. Finally, in Remark 4.13 we indicate how one can deal with exponential tail estimates for the weak limit.

We start by a lemma that ensures integrability properties needed in the rest of the paper (sometimes implicitly).

**Lemma 4.1.** *For a Borel function $g : \mathbb{R} \to \mathbb{R}$ such that $g(X_1) \in L_p$ with $p \in [2, \infty)$ and for $k, l \geq 0$, $j \in \{1, 2\}$, and $b \in [0, 1)$ one has that*

$$\mathbb{E} \sup_{0 \leq s \leq t \leq b} |X_t|^k |X_s|^l \left| \frac{\partial^j G}{\partial x^j}(s, X_s) \right|^p < \infty,$$

*where $G$ is given by* (2) *and* $0^0 := 1$.

The lemma is standard and proved in [7], Lemma 2.3, for $p = 2$ and integers $k, l \geq 0$ using the hyper-contractivity of the Ornstein–Uhlenbeck semigroup. Exactly the same proof works in our setting.

As before, we let $E = \mathbb{R}$ if $X = B$ and $E = (0, \infty)$ if $X = S$. Given a differentiable function $h : E \to \mathbb{R}$ we let

$$(Ah)(x) := (\sigma h')(x) - (\sigma' h)(x).$$

The main term is the first one; the second one guarantees that the operator $(Ag)(x)$ is constant in $x$ in the case $g(x) = c_0 x + c_1$ with $c_0, c_1 \in \mathbb{R}$, where the error process $C_t(\tau)$ of our approximation problem vanishes a.s. In the following $AG(t, x)$ always means that $A$ acts on the $x$-variable of the function $G(t, x)$.

**Definition 4.2.** *For $g(X_1) \in L_2$, $\beta \in (0, 1)$ and $t \in [0, 1)$ we let*

$$D_t^{X,\beta} g(X_1) := \frac{1-\beta}{2} \int_0^1 (1-u)^{-(1+\beta)/2} [AG(u \wedge t, X_{u \wedge t}) - AG(0, X_0)] \, du.$$

*For $\beta = 1$ and $t \in [0, 1)$ we let $D_t^{X,1} g(X_1) := AG(t, X_t) - AG(0, X_0)$.*

The process $D^{X,\beta}g(X_1) = (D_t^\beta g(X_1))_{t\in[0,1)}$ is a square-integrable martingale on the half open time interval $[0,1)$, because $((\sigma\frac{\partial G}{\partial x})(u, X_u))_{u\in[0,1)}$ and $(G(u, X_u))_{u\in[0,1]}$ are square-integrable martingales (cf. the remarks in the proof of Proposition 5.1). How should we interpret the case $\beta \in (0,1)$? Using the Riemann–Liouville operator of partial integration

$$(R_\alpha h)(t) := \frac{1}{\Gamma(\alpha)} \int_0^t (t-u)^{\alpha-1} h(u)\,du$$

for (say) continuous $h: [0,1] \to \mathbb{R}$ and $\alpha > 0$ we have

$$D_t^{X,\beta}g(X_1) = \Gamma\left(\frac{3-\beta}{2}\right)(R_{(1-\beta)/2}[AG(\cdot \wedge t, X_{\cdot \wedge t}) - AG(0, X_0)])(1).$$

That means that we differentiate once in the state direction by $A$ and integrate pathwise "back" $(1-\beta)/2$ times in time. Having in mind the parabolic PDE (3) this can be interpreted as integration in $x$ by an order $1-\beta$, so that we are left with a fractional differentiation of order $\beta$ in $x$. The point of the construction of $D^{X,\beta}g(X_1)$ is that we may have $L_p$-singularities of $(\sigma\frac{\partial G}{\partial x})(t, X_t)$ as $t \uparrow 1$ whereas $D^{X,\beta}g(X_1)$ stays $L_p$-bounded (see Example 4.7).

A first consequence of Theorem 3.1 is:

**Corollary 4.3.** *For $p \in [2, \infty)$, $\beta \in (0, 1]$ and $g(X_1) \in L_2$ the following assertions are equivalent:*

  (i) *On some stochastic basis there exists a continuous $L_p$-integrable martingale $M$ such that $\sqrt{n}C(\tau^{n,\beta}) \Longrightarrow_{C[0,1]} M$.*
  (ii) *The martingale $D^{X,\beta}g(X_1)$ is bounded in $L_p$.*

**Proof.** Let $\beta \in (0,1)$ and $t \in [0,1)$. By Itô's formula we get that

$$(1-t)^{(1-\beta)/2}\left[\left(\sigma\frac{\partial G}{\partial x}\right)(t, X_t) - (\sigma'G)(t, X_t)\right]$$
$$= \left[\left(\sigma\frac{\partial G}{\partial x}\right)(0, X_0) - (\sigma'G)(0, X_0)\right]$$
$$+ \int_0^t (1-u)^{(1-\beta)/2}\left(\sigma^2 \frac{\partial^2 G}{\partial x^2}\right)(u, X_u)\,dB_u$$
$$- \frac{1-\beta}{2}\int_0^t (1-u)^{-(1+\beta)/2}\left[\left(\sigma\frac{\partial G}{\partial x}\right)(u, X_u) - (\sigma'G)(u, X_u)\right]du.$$

By rearranging we arrive at

$$\int_0^t (1-u)^{(1-\beta)/2}\left(\sigma^2\frac{\partial^2 G}{\partial x^2}\right)(u, X_u)\,dB_u$$

$$= -\left[\left(\sigma\frac{\partial G}{\partial x}\right)(0, X_0) - (\sigma' G)(0, X_0)\right]$$

$$+ \frac{1-\beta}{2}\int_0^1 (1-u)^{-(1+\beta)/2}\left[\left(\sigma\frac{\partial G}{\partial x}\right)(u\wedge t, X_{u\wedge t}) - (\sigma' G)(u\wedge t, X_{u\wedge t})\right]du.$$

Applying the Burkholder–Davis–Gundy inequalities we deduce that

$$\left(\int_0^1 (1-u)^{1-\beta}\left[\left(\sigma^2\frac{\partial^2 G}{\partial x^2}\right)(u, X_u)\right]^2 du\right)^{1/2} \in L_p \tag{9}$$

if and only if Assertion (ii) of our theorem is satisfied. The equivalence of (9) to (i) follows from Theorem 3.1, (5) and $\mathbb{E}|W_A|^p = \mathbb{E}|A|^{p/2}\mathbb{E}|W_1|^p$ for $A := A_\beta(1)\chi_{\{A_\beta(1)<\infty\}}$. The case $\beta = 1$ can be treated in a similar way. $\square$

For $\beta = 1$ it follows from Corollary 4.3 that $\sup_{t\in[0,1)} \|AG(t, X_t)\|_{L_p} < \infty$ is equivalent to the existence of a continuous $L_p$-integrable martingale $M$ such that $\sqrt{n}C(\tau^{n,1}) \Longrightarrow_{C[0,1]} M$. Now we treat the case $\beta \in (0, 1)$.

**Corollary 4.4.** *For $g(X_1) \in L_2$ and $p \in [2, \infty)$ one has the following:*

(i) *If $0 < \beta < \alpha \leq 1$ and*

$$\sup_{t\in[0,1)} (1-t)^{(1-\alpha)/2}\|AG(t, X_t)\|_{L_p} < \infty,$$

*then on some stochastic basis there exists a continuous $L_p$-integrable martingale $M$ such that $\sqrt{n}C(\tau^{n,\beta}) \Longrightarrow_{C[0,1]} M$.*

(ii) *If $\beta \in (0, 1]$ and if there exists a continuous $L_p$-integrable martingale $M$ such that $\sqrt{n}C(\tau^{n,\beta}) \Longrightarrow_{C[0,1]} M$, then*

$$\sup_{t\in[0,1)} (1-t)^{(1-\beta)/2}\|AG(t, X_t)\|_{L_p} < \infty.$$

**Proof.** (i) There exists a $c > 0$ such that, for $t \in [0, 1)$, we have that

$$\|D_t^{X,\beta}g(X_1)\|_{L_p}$$

$$= \frac{1-\beta}{2}\left\|\int_0^1 (1-u)^{-(1+\beta)/2}[AG(u\wedge t, X_{u\wedge t}) - AG(0, X_0)]du\right\|_{L_p}$$

$$\leq \frac{1-\beta}{2}\int_0^1 (1-u)^{-(1+\beta)/2}\|AG(u, X_u)\|_{L_p}du + |AG(0, X_0)|$$

$$\leq c\frac{1-\beta}{2}\int_0^1 (1-u)^{-(1+\beta)/2}(1-u)^{(\alpha-1)/2}du + |AG(0, X_0)|$$

$$< \infty$$

and can apply Corollary 4.3.

(ii) From Theorem 3.1 we get that

$$\left\| \left( \int_0^1 (1-u)^{1-\beta} \left[ \left( \sigma^2 \frac{\partial^2 G}{\partial x^2} \right)(u, X_u) \right]^2 du \right)^{1/2} \right\|_{L_p} < \infty$$

and, using the Burkholder–Davis–Gundy inequalities, the existence of a constant $c > 0$ such that

$$\left\| \int_0^t \left( \sigma^2 \frac{\partial^2 G}{\partial x^2} \right)(u, X_u) \, dB_u \right\|_{L_p} \le c(1-t)^{(\beta-1)/2}$$

for $t \in [0, 1)$. But the left-hand side can be re-written as $AG(t, X_t) - AG(0, X_0)$ and we are done. □

Now we extend [6], Lemma 3.6, from $p = 2$ to $p \in (2, \infty)$.

**Proposition 4.5.** *For $p \in [2, \infty)$, $\beta \in (0, 1]$ and $g(X_1) \in L_p$ the following assertions are equivalent:*

(i) $\sup_{t \in [0,1)} (1-t)^{(1-\beta)/2} \|AG(t, X_t)\|_{L_p} < \infty$.
(ii) $\sup_{t \in [0,1)} (1-t)^{-\beta/2} \|g(X_1) - \mathbb{E}(g(X_1)|\mathcal{F}_t)\|_{L_p} < \infty$.
(iii) *There exists a constant $c > 0$ such that for all two-dimensional Gaussian random vectors $(Y, Z)$ with $Y, Z \sim N(0, 1)$ one has*

$$\mathbb{E}|h(Y) - h(Z)|^p \le c\mathbb{E}|Y - Z|^{\beta p},$$

*where $h := g$ if $X = B$ and $h(x) := g(e^{x-1/2})$ if $X = S$.*

From (iii) it is clear that Hölder continuous functions $h$ with exponent $\beta$ satisfy the properties of Proposition 4.5. But, for example, one also has

$$\mathbb{E}|\chi_{[K,\infty)}(Y) - \chi_{[K,\infty)}(Z)|^p \le c\mathbb{E}|Y - Z|^{\beta p}$$

for $K \in \mathbb{R}$ and $\beta = \frac{1}{p}$, where $c > 0$ is independent from $Y$ and $Z$, as shown in Example 4.7 below.

**Proof of Proposition 4.5.** By $g(X_1) \in L_p$ it follows that assertion (i) is equivalent to the existence of a $c > 0$ such that

$$\left\| \left( \sigma \frac{\partial G}{\partial x} \right)(t, X_t) \right\|_{L_p} \le c(1-t)^{(\beta-1)/2}.$$

(i) $\implies$ (ii) is clear as $p \in [2, \infty)$ an interchange of the $L_p$- and $L_2$-norm and the Burkholder–Davis–Gundy inequalities give that

$$\|g(X_1) - \mathbb{E}(g(X_1)|\mathcal{F}_t)\|_{L_p} \leq c_p \left( \int_t^1 \left\| \left( \sigma \frac{\partial G}{\partial x} \right)(s, X_s) \right\|_{L_p}^2 ds \right)^{1/2}$$

$$\leq c_p c \left( \int_t^1 (1-s)^{\beta-1} ds \right)^{1/2}$$

$$\leq \frac{c_p c}{\sqrt{\beta}} (1-t)^{\beta/2}.$$

(ii) $\implies$ (i) Here it is known (see [19] and [9], Lemmas A.1 and A.2) that for $t \in [0, 1)$, a.s.,

$$\left( \sigma \frac{\partial G}{\partial x} \right)(t, X_t) = \mathbb{E}\left( g(X_1) \frac{B_1 - B_t}{1 - t} \Big| \mathcal{F}_t \right)$$

$$= \mathbb{E}\left( [g(X_1) - \mathbb{E}(g(X_1)|\mathcal{F}_t)] \frac{B_1 - B_t}{1 - t} \Big| \mathcal{F}_t \right).$$

For $1 = \frac{1}{p} + \frac{1}{q}$ and $c_q := \|B_1\|_{L_q}$ this implies that, a.s.,

$$\left| \left( \sigma \frac{\partial G}{\partial x} \right)(t, X_t) \right|$$

$$\leq (\mathbb{E}(|g(X_1) - \mathbb{E}(g(X_1)|\mathcal{F}_t)|^p | \mathcal{F}_t))^{1/p} \left( \mathbb{E}\left( \left| \frac{B_1 - B_t}{1 - t} \right|^q \Big| \mathcal{F}_t \right) \right)^{1/q}$$

$$= (\mathbb{E}(|g(X_1) - \mathbb{E}(g(X_1)|\mathcal{F}_t)|^p | \mathcal{F}_t))^{1/p} c_q (1 - t)^{-1/2}.$$

By integration the desired inequality follows since

$$\left\| \left( \sigma \frac{\partial G}{\partial x} \right)(t, X_t) \right\|_{L_p} \leq c_q (1-t)^{-1/2} \|g(X_1) - \mathbb{E}(g(X_1)|\mathcal{F}_t)\|_{L_p}$$

$$\leq c_q c (1-t)^{(\beta-1)/2}.$$

(iii) is a reformulation of (ii) because

$$\|g(X_1) - G(t, X_t)\|_{L_p}$$
$$= \|h(B_1) - \mathbb{E}(h(B_1)|\mathcal{F}_t)\|_{L_p}$$
$$= \|h(B_1) - \widetilde{\mathbb{E}} h(B_t + \widetilde{B}_{1-t})\|_{L_p}$$
$$\leq \|h(B_1) - h(B_t + \widetilde{B}_{1-t})\|_{L_p}$$
$$\leq \|h(B_1) - \mathbb{E}(h(B_1)|\mathcal{F}_t)\|_{L_p} + \|\mathbb{E}(h(B_1)|\mathcal{F}_t) - h(B_t + \widetilde{B}_{1-t})\|_{L_p}$$

$$= 2\|g(X_1) - G(t, X_t)\|_{L_p},$$

where $\widetilde{B}$ is a Brownian motion on an auxiliary stochastic basis with expected value $\widetilde{E}$. For the computation above we use that $(B_1, B_t + \widetilde{B}_{1-t})$ is distributed like a two-dimensional Gaussian random vector $(Y, Z)$ with $Y, Z \sim N(0,1)$ and $\mathrm{cov}(Y, Z) = t$, and that $\|Y - Z\|_{L_{\beta p}} \sim \sqrt{1 - \mathrm{cov}(Y, Z)}$ (the same argument was exploited in [12], proof of Corollary 2.3). $\square$

**Remark 4.6.** Proposition 4.5(iii) gives the direct link to the spaces $\mathcal{E}_p^\alpha$ considered in [16], page 428, with

$$\mathcal{E}_p^\beta = \left\{ g \in L_p(\gamma) : \int_0^\infty t^{-1-\beta p/2} \mathbb{E}|g(g_1) - g(\mathrm{e}^{-t/2}g_1 + \sqrt{1-\mathrm{e}^{-t}}g_2)|^p \, \mathrm{d}t < \infty \right\},$$

where $\beta \in (0,1)$, $p \in [2, \infty)$, and $g_1$ and $g_2$ are independent standard Gaussian random variables. There is a slight difference between these spaces and the condition we use: for $X = B$ and $g(x) = \chi_{[K,\infty)}(x)$ we have that

$$\left\| \frac{\partial G}{\partial x}(t, B_t) \right\|_{L_p} \sim (1-t)^{1/(2p) - 1/2}$$

as shown in Example 4.7 below. This implies that the conditions of Proposition 4.5 are satisfied for $\beta = \frac{1}{p}$. The arguments of the proof of Proposition 4.5 also give

$$\mathbb{E}|g(g_1) - g(\mathrm{e}^{-t/2}g_1 + \sqrt{1-\mathrm{e}^{-t}}g_2)|^p \sim (1-\mathrm{e}^{-t/2})^{1/2},$$

so that

$$\chi_{[K,\infty)} \notin \mathcal{E}_p^{1/p} \quad \text{but} \quad \chi_{[K,\infty)} \in \mathcal{E}_p^\beta \qquad \text{whenever } 0 < \beta < \frac{1}{p}.$$

**Example 4.7.** Let $g(x) := \chi_{[K,\infty)}(x)$ where $K > 0$ if $X = S$ and let $\beta \in (0, \frac{1}{2})$. Then one has

$$\sup\left\{ p \in [2, \infty) : \text{the weak limit } \lim_n \sqrt{n} C(\tau^{n,\beta}) \text{ is } L_p\text{-integrable} \right\} = \frac{1}{\beta}.$$

**Proof.** We have only to consider the case $X = B$: assuming $g_B(B_1) = \chi_{[K_B,\infty)}(B_1) = \chi_{[K_S,\infty)}(S_1) = g_S(S_1)$, one gets for the solutions of the corresponding backward PDEs, $G_B$ and $G_S$, that $\frac{\partial G_B}{\partial x}(t, B_t) = S_t \frac{\partial G_S}{\partial x}(t, S_t)$. Hence, $D^{B,\beta} g_B(B_1)$ is bounded in $L_p$ if and only if $D^{S,\beta} g_S(S_1)$ is bounded in $L_p$. Therefore we may assume that $X = B$ so that

$$G(t, x) = \int_{K-x}^\infty \mathrm{e}^{-y^2/(2(1-t))} \frac{\mathrm{d}y}{\sqrt{2\pi(1-t)}} \quad \text{and} \quad \frac{\partial G}{\partial x}(t, x) = \mathrm{e}^{-(K-x)^2/(2(1-t))} \frac{1}{\sqrt{2\pi(1-t)}}$$

for $t \in [0,1)$. Then for $p \in [2,\infty)$ and $t \in [0,1)$ we get that

$$\left\|\frac{\partial G}{\partial x}(t, B_t)\right\|_{L_p} \sim (1-t)^{1/(2p)-1/2},$$

which gives the assertion by Corollary 4.4. $\square$

**Example 4.8.** If $g(X_1) \in \bigcap_{p \in [2,\infty)} L_p$ and if there exist $\theta \in [0, 1/2)$ and $q, c \in (0, \infty)$ such that

$$\left|\left(\sigma \frac{\partial G}{\partial x}\right)(t, x)\right| \leq c \frac{1+|x|^q}{(1-t)^\theta},$$

then

$$\sup\left\{p \in [2,\infty) : \text{the weak limit } \lim_n \sqrt{n} C(\tau^{n,\beta}) \text{ is } L_p\text{-integrable}\right\} = \infty$$

for $\beta \in (0, 1-2\theta)$. In particular, it holds for $X = S$, Hölder continuous $g:\mathbb{R} \to \mathbb{R}$ of exponent $\eta \in (0,1]$ and $\theta = \frac{1-\eta}{2}$.

**Proof.** The first part follows by Corollary 4.4 with $\alpha := 1 - 2\theta$. So assume $|g(x) - g(y)| \leq d|x-y|^\eta$ for some $d > 0$ and $X = S$. Then it is known that

$$\left|x\frac{\partial G}{\partial x}(t, x)\right| = \left|\mathbb{E}g(xS_{1-t})\frac{B_{1-t}}{1-t}\right| = \left|\mathbb{E}[g(xS_{1-t}) - g(x)]\frac{B_{1-t}}{1-t}\right|$$
$$\leq d'|x|^\eta (\mathbb{E}|S_{1-t}-1|^2)^{\eta/2}(1-t)^{-1/2} \leq d''|x|^\eta (1-t)^{(\eta-1)/2}. \quad \square$$

Finally, we exploit Proposition 4.5(iii) to show that local properties of $g$ yield upper bounds for the integrability of the weak limit from Theorem 3.1. The local properties of $g$ are formulated by the following version of a sharp function: given $p \in [1,\infty)$, a locally $L_p$-integrable $g: \mathbb{R} \to \mathbb{R}$, $x_0 \in \mathbb{R}$ and $\varepsilon > 0$, we let

$$\mathrm{OSC}_p(g, x_0, \varepsilon) := \left(\frac{1}{4\varepsilon^2}\int_{Q(x_0,\varepsilon)} |g(y)-g(z)|^p \, \mathrm{d}y\,\mathrm{d}z\right)^{1/p},$$

where $Q(x_0, \varepsilon) := \{(y,z) : |y-x_0| \leq \varepsilon, |z-x_0| \leq \varepsilon\} \subseteq \mathbb{R}^2$.

**Lemma 4.9.** For all $p \in [1,\infty)$ and $x_0 \in \mathbb{R}$ there is a constant $c > 0$ such that

$$\mathrm{OSC}_p(g, x_0, \sqrt{1-t}) \leq c(1-t)^{-1/(2p)}\|g(Y) - g(Z)\|_{L_p}$$

for all $t \in [0,1)$, $g \in L_p(\gamma)$, and two-dimensional Gaussian random vectors $(Y, Z)$ with $Y, Z \sim N(0,1)$ and $\mathrm{cov}(Y, Z) = t$.

**Proof.** Given $t \in [0, 1)$, let $(Y, Z)$ be the above two-dimensional Gaussian random vector, that is, we have the covariance

$$C_t := \begin{pmatrix} 1 & t \\ t & 1 \end{pmatrix} \quad \text{so that} \quad C_t^{-1} := \frac{1}{1-t^2} \begin{pmatrix} 1 & -t \\ -t & 1 \end{pmatrix}$$

and the density of the law of $(Y, Z)$ can be computed as

$$p_t(y, z) := \frac{1}{2\pi\sqrt{1-t^2}} e^{-(y^2 - 2tyz + z^2)/(2(1-t^2))}.$$

Letting $(y, z) \in Q(x_0, \sqrt{1-t})$ with $y = x_0 + r$, $z = x_0 + s$ and $r, s \in [-\sqrt{1-t}, \sqrt{1-t}]$ we get that

$$y^2 - 2tyz + z^2 = 2(1-t)x_0[x_0 + s + r] + r^2 + s^2 - 2trs$$

so that $|y^2 - 2tyz + z^2| \le A(1-t)$ with $A = A(x_0) > 0$ and

$$p_t(y, z) \ge \frac{1}{2\pi\sqrt{1-t^2}} e^{-A/(2(1+t))} \ge \frac{1}{16 e^A \sqrt{1-t}} \quad \text{for } (y, z) \in Q(x_0, \sqrt{1-t}).$$

Consequently,

$$\mathbb{E}|g(Y) - g(Z)|^p \ge \frac{1}{16 e^A \sqrt{1-t}} \int_{Q(x_0, \sqrt{1-t})} |g(y) - g(z)|^p \, dy \, dz$$

$$= \frac{1}{4 e^A} \sqrt{1-t} \left( \frac{1}{4(1-t)} \int_{Q(x_0, \sqrt{1-t})} |g(y) - g(z)|^p \, dy \, dz \right). \qquad \square$$

*Remark 4.10.* We use Lemma 4.9 implicitly in the proof of Corollary 4.11 below to get a statement (in a sense) opposite to [25], Corollary on page 185: Instead of deducing Hölder continuity properties of $g$ we assume that the local oscillation of $g$ in $x_0 \in \mathbb{R}$ is singular of order $\eta \ge 0$ to deduce the upper bound $\beta \le \frac{1}{p} - \eta$ for our fractional smoothness.

**Corollary 4.11.** *Let $X = B$, $p \in [2, \infty)$, $\eta \in [0, \frac{1}{2})$, $\beta \in (0, 1]$ and $g \in L_p(\gamma)$.*

(i) *If there exists a continuous $L_p$-integrable martingale $M$ such that*

$$\sqrt{n} C(\tau^{n,\beta}) \Longrightarrow_{C[0,1]} M$$

*and*

(ii) *if there exists an $x_0 \in \mathbb{R}$ such that $\limsup_{1 \ge \varepsilon \to 0} \varepsilon^\eta \operatorname{OSC}_p(g, x_0, \varepsilon) > 0$,*

*then, necessarily, $p \le \frac{1}{\beta + \eta}$.*

Note that $\operatorname{OSC}_p(g, x_0, \varepsilon)$ is monotone in $p$, so that we can use $\operatorname{OSC}_2(g, x_0, \varepsilon)$ in (ii) as well.

**Proof of Corollary 4.11.** By Corollary 4.4, Assumption (i) ensures that

$$\sup_{t\in[0,1)} (1-t)^{(1-\beta)/2}\|AG(t,B_t)\|_{L_p} < \infty,$$

so that the statements of Proposition 4.5 are valid for this $\beta$. On the other hand, Assumption (ii) ensures a constant $d > 0$ and a sequence $t_n \uparrow 1$ such that, by Lemma 4.9 and Proposition 4.5(iii),

$$\frac{1}{d}(\sqrt{1-t_n})^{-\eta} \le \mathrm{OSC}_p(g, x_0, \sqrt{1-t_n})$$
$$\le c_{(4.9)}(1-t_n)^{-1/(2p)}\|g(Y_n) - g(Z_n)\|_{L_p}$$
$$\le c_{(4.9)}c_{(4.5)}(1-t_n)^{-1/(2p)}(1-t_n)^{\beta/2},$$

where $\mathrm{cov}(Y_n, Z_n) = t_n$. But this implies $\beta + \eta \le \frac{1}{p}$. □

***Example 4.12.*** (i) *Let $p \in [1, \infty)$, $\eta \in [0, 1/p)$ and define $g(x) := 0$ for $x \le 0$ and $g(x) := x^{-\eta}$ for $x > 0$. Then*

$$\frac{1}{c}\varepsilon^{-\eta} \le \mathrm{OSC}_p(g, 0, \varepsilon) \le c\varepsilon^{-\eta} \qquad \textit{for } \varepsilon \in (0,1] \textit{ where } c = c(p,\eta) > 0.$$

(ii) *Let $g : \mathbb{R} \to \mathbb{R}$ be locally $L_p$-integrable, $p \in [1, \infty)$, such that there are $x_0, y_0 \in \mathbb{R}$ and $\delta > 0$ with $g(x) \le y_0 - \delta$ for $x \in (x_0 - \delta, x_0)$ and $g(x) \ge y_0 + \delta$ for $x \in (x_0, x_0 + \delta)$. Then*

$$\limsup_{1 \ge \varepsilon \to 0} \mathrm{OSC}_p(g, x_0, \varepsilon) \ge 2^{1-1/p}\delta.$$

***Remark 4.13.*** So far we considered $L_p$-bounds for the weak limit and used the Burkholder–Davis–Gundy inequalities that preserve these $L_p$-bounds. In the limiting case $p = \infty$ we have to proceed differently and exploit for a random variable $Z$ the following equivalence: given $r \in (0, \infty)$ there is a $c > 0$ such that

$$\mathbb{P}(|Z| > \lambda) \le ce^{-\lambda^r/c} \qquad \text{for } \lambda > 0 \quad \text{if and only if} \quad \sup_{p\in[1,\infty)} p^{-1/r}\|Z\|_{L_p} < \infty.$$

Letting $A := A_\beta(1)\chi_{\{A_\beta(1)<\infty\}}$ and applying this equivalence to

$$\|W_A\|_{L_p} = \|\sqrt{A}\|_{L_p}\|W_1\|_{L_p} \sim \sqrt{p}\|\sqrt{A}\|_{L_p}$$

for $p \in [1, \infty)$, one can compare the tail behavior of $W_A$ and that of $\sqrt{A}$. In particular,

$$\|A\|_{L_\infty} < \infty \tag{10}$$

if and only if there exists a constant $c > 0$ such that

$$\mathbb{P}(|W_{A_\beta(1)\chi_{\{A_\beta(1)<\infty\}}}| > \lambda) \leq c e^{-\lambda^2/c} \qquad \text{for } \lambda > 0.$$

A typical example for (10) is obtained if there exist $\theta \in [0,1)$ and $c > 0$ such that

$$\left|\left(\sigma^2 \frac{\partial^2 G}{\partial x^2}\right)(t,x)\right| \leq \frac{c}{(1-t)^\theta}. \tag{11}$$

Then obviously $\|A_\beta(1)\|_{L_\infty} < \infty$ is satisfied for $\beta \in (0, 2(1-\theta)) \cap (0,1]$. In particular, condition (11) holds for $X = B$ and $\eta$-Hölder continuous $g : \mathbb{R} \to \mathbb{R}$ and for $X = S$ if $g(e^x)$ is $\eta$-Hölder continuous, where $\theta = 1 - \frac{\eta}{2} \in [\frac{1}{2}, 1)$, $\eta \in (0,1]$. Here one can follow the same pattern as in Example 4.8 because

$$\frac{\partial^2 G}{\partial x^2}(t,x) = \mathbb{E} g(x + \sqrt{1-t} B_1) \frac{B_1^2 - 1}{1-t}$$

for $X = B$ and

$$x^2 \frac{\partial^2 G}{\partial x^2}(t,x) = \mathbb{E}\left(g(xS_{1-t})\left[\frac{B_{1-t}^2 - (1-t)}{(1-t)^2} - \frac{B_{1-t}}{1-t}\right]\right)$$

for $X = S$.

## 5. Proof of Theorem 3.1

Throughout this section we let

$$H(t) := \left\|\left(\sigma^2 \frac{\partial^2 G}{\partial x^2}\right)(t, X_t)\right\|_{L_2} \qquad \text{for } t \in [0,1)$$

and obtain a continuous and non-decreasing function $H : [0,1) \to [0,\infty)$ (see [12], Lemma 3.9).

### 5.1. Proof of (iii) $\Longrightarrow$ (ii) $\Longrightarrow$ (i)

For the first implication we simply take $M_t := \widetilde{Z}_\beta(t)$ and the natural filtration of $M$. For the second one we remark that (ii) implies that

$$\frac{1}{2}\int_0^1 \nu_\beta(u) H(u)^2 \, du = \sup_{t \in [0,1)} \mathbb{E} \frac{1}{2} \int_0^t \nu_\beta(u) \left[\left(\sigma^2 \frac{\partial^2 G}{\partial x^2}\right)(u, X_u)\right]^2 du$$

$$= \sup_{t \in [0,1)} \mathbb{E} Z_\beta(t)^2$$

$$= \sup_{t \in [0,1)} \mathbb{E} M_t^2$$

$$\leq \mathbb{E} M_1^2 < \infty$$

so that

$$\int_0^1 (1-u)^{1-\beta} H(u)^2 \, du < \infty.$$

Now Assertion (i) follows from (5).

## 5.2. Preparations for the proof of (i) $\implies$ (iii)

First we decompose the error process. For $t \in [0,1]$ and a time-net $\tau = (t_i)_{i=0}^n$, $0 = t_0 < \cdots < t_n = 1$, we obtain, $\mathbb{P}$-a.s., that

$$\begin{aligned}
C_t(\tau) &= \left[ \int_0^t \frac{\partial G}{\partial x}(u, X_u) \, dX_u - \sum_{i=0}^{n-1} \frac{\partial G}{\partial x}(t_i, X_{t_i})(X_{t_{i+1} \wedge t} - X_{t_i \wedge t}) \right] \\
&= \sum_{i=0}^{n-1} \int_{t_i \wedge t}^{t_{i+1} \wedge t} \left[ \frac{\partial G}{\partial x}(u, X_u) - \frac{\partial G}{\partial x}(t_i, X_{t_i}) \right. \\
&\qquad \left. - \frac{\partial^2 G}{\partial x^2}(t_i, X_{t_i})(X_u - X_{t_i}) \right] dX_u \\
&\quad + \sum_{i=0}^{n-1} \int_{t_i \wedge t}^{t_{i+1} \wedge t} [\sigma(X_u) - \sigma(X_{t_i})] \frac{\partial^2 G}{\partial x^2}(t_i, X_{t_i})(X_u - X_{t_i}) \, dB_u \\
&\quad + \sum_{i=0}^{n-1} \int_{t_i \wedge t}^{t_{i+1} \wedge t} \left( \sigma \frac{\partial^2 G}{\partial x^2} \right)(t_i, X_{t_i})(X_u - X_{t_i}) \, dB_u \\
&=: I_t^1(\tau) + I_t^2(\tau) + I_t^3(\tau).
\end{aligned}$$

The appropriate $L_2$-integrability of the integrands in the decomposition above is obtained by standard arguments (see, e.g., Lemma 4.1 and its proof).

*Estimation of $I^1(\tau)$ and $I^2(\tau)$.*

First we show that a Taylor expansion of order one of the integrand of the stochastic integral $\int_0^1 \frac{\partial G}{\partial x}(u, X_u) \, dX_u$ gives an $L_2$-approximation rate of $o(1/\sqrt{n})$ provided that appropriate time-nets are taken.

**Proposition 5.1.** *Let $\beta \in (0,1]$, and $g \in B_{2,2}^{\beta}(\gamma)$ for $X = B$ and $g(e^{\cdot -1/2}) \in B_{2,2}^{\beta}(\gamma)$ for $X = S$, respectively. Then one has that*

$$\lim_n n\mathbb{E}|I_1^1(\tau^{n,\beta})|^2 = \lim_n n\mathbb{E}\left|\sum_{i=0}^{n-1}\int_{t_i^{n,\beta}}^{t_{i+1}^{n,\beta}}\left[\frac{\partial G}{\partial x}(u, X_u) - \frac{\partial G}{\partial x}(t_i^{n,\beta}, X_{t_i^{n,\beta}})\right.\right.$$
$$\left.\left. - \frac{\partial^2 G}{\partial x^2}(t_i^{n,\beta}, X_{t_i^{n,\beta}})(X_u - X_{t_i^{n,\beta}})\right]\mathrm{d}X_u\right|^2$$
$$= 0.$$

**Proof.** (a) Let $0 \le a < b < 1$ and

$$\Phi(u,x) = \Phi(u,x,\omega)$$
$$:= \left[\frac{\partial G}{\partial x}(u,x) - \frac{\partial G}{\partial x}(a, X_a(\omega)) - \frac{\partial^2 G}{\partial x^2}(a, X_a(\omega))(x - X_a(\omega))\right]\sigma(x)$$

for $u \in [a,b]$ and $x \in E$, where we shall suppress $\omega$ in the following. By a computation we get that

$$\frac{\partial \Phi}{\partial u}(u,x) + \frac{\sigma(x)^2}{2}\frac{\partial^2 \Phi}{\partial x^2}(u,x) = -\frac{\partial^2 G}{\partial x^2}(a, X_a)\sigma'(x)\sigma(x)^2.$$

From this we conclude that

$$\frac{\partial \Phi^2}{\partial u}(u,x) + \frac{\sigma(x)^2}{2}\frac{\partial^2 \Phi^2}{\partial x^2}(u,x)$$
$$= 2\Phi(u,x)\left[\frac{\partial \Phi}{\partial u}(u,x) + \frac{\sigma(x)^2}{2}\frac{\partial^2 \Phi}{\partial x^2}(u,x)\right] + \sigma(x)^2\left[\frac{\partial \Phi}{\partial x}(u,x)\right]^2$$
$$= -2\Phi(u,x)\frac{\partial^2 G}{\partial x^2}(a, X_a)\sigma'(x)\sigma(x)^2 + \sigma(x)^2\left[\frac{\partial \Phi}{\partial x}(u,x)\right]^2$$
$$= -2\Phi(u,x)\frac{\partial^2 G}{\partial x^2}(a, X_a)\sigma'(x)\sigma(x)^2$$
$$+ \left[\Phi(u,x)\sigma'(x) + \left[\frac{\partial^2 G}{\partial x^2}(u,x) - \frac{\partial^2 G}{\partial x^2}(a, X_a)\right]\sigma(x)^2\right]^2.$$

This yields

$$\left|\frac{\partial \Phi^2}{\partial u}(u,x) + \frac{\sigma(x)^2}{2}\frac{\partial^2 \Phi^2}{\partial x^2}(u,x)\right| \le 2\left|\Phi(u,x)\frac{\partial^2 G}{\partial x^2}(a, X_a)\sigma(x)^2\right|$$
$$+ 2\Phi(u,x)^2 + 2\left[\frac{\partial^2 G}{\partial x^2}(u,x) - \frac{\partial^2 G}{\partial x^2}(a, X_a)\right]^2\sigma(x)^4.$$

Moreover, by Itô's formula,

$$\mathbb{E}\Phi(b,X_b)^2 = \mathbb{E}\Phi(a,X_a)^2 + \mathbb{E}\int_a^b \left[\frac{\partial \Phi^2}{\partial u}(u,X_u) + \frac{\sigma(X_u)^2}{2}\frac{\partial^2 \Phi^2}{\partial x^2}(u,X_u)\right]du.$$

In fact, first we condition on $X_a = y$, then we apply Itô's formula on $[a,b]$ to obtain (conditionally) the equation with $b$ replaced by $\tau_N$ defined as the minimum of $b$, $\inf\{s \in [a,b] : |(\partial(\Phi^2)/\partial x)(u,X_u)| \geq N\}$ and $\inf\{s \in [a,b] : |X_u - y| \geq N$ if $X = B$, $|X_u/y| \notin ((1/N), N)$ if $X = S\}$. Finally we let $N \to \infty$ by the help of Lemma 4.1 (cf. [10], proof of Theorem 6, for the conditioning argument). One more integration to remove the condition $X_a = y$ gives the equation we want.

Now Gronwall's lemma gives

$$\mathbb{E}\Phi(b,X_b)^2 \leq c_{(12)}\left[\int_a^b \mathbb{E}\left|\Phi(u,X_u)\frac{\partial^2 G}{\partial x^2}(a,X_a)\sigma(X_u)^2\right|du \right. \\ \left. + \int_a^b \mathbb{E}\left[\frac{\partial^2 G}{\partial x^2}(u,X_u) - \frac{\partial^2 G}{\partial x^2}(a,X_a)\right]^2 \sigma(X_u)^4\,du\right] \quad (12)$$

for some absolute constant $c_{(12)} > 0$.

(b) Let $i \in \{0,\ldots,n-1\}$ and $u \in [t_i^{n,\beta}, t_{i+1}^{n,\beta})$ and set

$$\Phi_i^n(u,x) := \left[\frac{\partial G}{\partial x}(u,x) - \frac{\partial G}{\partial x}(t_i^{n,\beta}, X_{t_i^{n,\beta}}) - \frac{\partial^2 G}{\partial x^2}(t_i^{n,\beta}, X_{t_i^{n,\beta}})(x - X_{t_i^{n,\beta}})\right]\sigma(x).$$

From step (a) we conclude that

$$n\mathbb{E}\left|\sum_{i=0}^{n-1}\int_{t_i^{n,\beta}}^{t_{i+1}^{n,\beta}}\left[\frac{\partial G}{\partial x}(u,X_u) - \frac{\partial G}{\partial x}(t_i^{n,\beta}, X_{t_i^{n,\beta}})\right.\right.$$

$$\left.\left. - \frac{\partial^2 G}{\partial x^2}(t_i^{n,\beta}, X_{t_i^{n,\beta}})(X_u - X_{t_i^{n,\beta}})\right]dX_u\right|^2$$

$$= n\sum_{i=0}^{n-1}\int_{t_i^{n,\beta}}^{t_{i+1}^{n,\beta}} \mathbb{E}\Phi_i^n(u,X_u)^2\,du$$

$$\leq c_{(12)}n\sum_{i=0}^{n-1}\int_{t_i^{n,\beta}}^{t_{i+1}^{n,\beta}} \left[\int_{t_i^{n,\beta}}^u \mathbb{E}\left|\Phi_i^n(v,X_v)\frac{\partial^2 G}{\partial x^2}(t_i^{n,\beta}, X_{t_i^{n,\beta}})\sigma(X_v)^2\right|dv\right.$$

$$\left. + \int_{t_i^{n,\beta}}^u \mathbb{E}\left[\frac{\partial^2 G}{\partial x^2}(v,X_v) - \frac{\partial^2 G}{\partial x^2}(t_i^{n,\beta}, X_{t_i^{n,\beta}})\right]^2 \sigma(X_v)^4\,dv\right]du$$

$$= c_{(12)}n\sum_{i=0}^{n-1}\int_{t_i^{n,\beta}}^{t_{i+1}^{n,\beta}}\int_{t_i^{n,\beta}}^u A_i^n(v)^2\,dv\,du$$

with

$$A_i^n(v)^2 := \mathbb{E}\left|\Phi_i^n(v,X_v)\frac{\partial^2 G}{\partial x^2}(t_i^{n,\beta},X_{t_i^{n,\beta}})\sigma(X_v)^2\right|$$
$$+ \mathbb{E}\left[\frac{\partial^2 G}{\partial x^2}(v,X_v) - \frac{\partial^2 G}{\partial x^2}(t_i^{n,\beta},X_{t_i^{n,\beta}})\right]^2 \sigma(X_v)^4$$

for $v \in [t_i^{n,\beta}, t_{i+1}^{n,\beta})$. Using (4) we continue by (cf. [12])

$$c_{(12)}n\sum_{i=0}^{n-1}\int_{t_i^{n,\beta}}^{t_{i+1}^{n,\beta}}\int_{t_i^{n,\beta}}^{u}A_i^n(v)^2\,dv\,du = c_{(12)}n\sum_{i=0}^{n-1}\int_{t_i^{n,\beta}}^{t_{i+1}^{n,\beta}}(t_{i+1}^{n,\beta}-u)A_i^n(u)^2\,du$$
$$\leq \frac{c_{(12)}}{\beta}\sum_{i=0}^{n-1}\int_{t_i^{n,\beta}}^{t_{i+1}^{n,\beta}}(1-u)^{1-\beta}A_i^n(u)^2\,du$$
$$= \frac{c_{(12)}}{\beta}\int_0^1 (1-u)^{1-\beta}\psi_n(u)\,du$$

with

$$\psi_n(u) := \sum_{i=0}^{n-1}\chi_{[t_i^{n,\beta},t_{i+1}^{n,\beta})}(u)A_i^n(u)^2.$$

(c) Now we show that

$$\psi_n(u) \leq c_{(13)}|H(u) \vee \|g(X_1)\|_{L_2}|^2 \tag{13}$$

for some absolute constant $c_{(13)} > 0$. Assume again $a = t_i^{n,\beta} \leq u < t_{i+1}^{n,\beta}$. Since the process $((\sigma^2 \frac{\partial^2 G}{\partial x^2})(u,X_u))_{u\in[0,1)} \subseteq L_2$ is a martingale (the argument for $X = S$ is given in [8]; the case $X = B$ can be treated in the same way) we get that

$$\mathbb{E}\left[\frac{\partial^2 G}{\partial x^2}(u,X_u) - \frac{\partial^2 G}{\partial x^2}(a,X_a)\right]^2 \sigma(X_u)^4 \leq c_{(14)}H(u)^2 \tag{14}$$

for some absolute constant $c_{(14)} > 0$. The first term of $A_i^n(u)^2$ can be bounded by

$$\mathbb{E}\left|\Phi_i^n(u,X_u)\sigma(X_u)^2\frac{\partial^2 G}{\partial x^2}(a,X_a)\right|$$
$$= \mathbb{E}\left|\left[\frac{\partial G}{\partial x}(u,X_u) - \frac{\partial G}{\partial x}(a,X_a) - \frac{\partial^2 G}{\partial x^2}(a,X_a)(X_u-X_a)\right]\sigma(X_u)\sigma(X_u)^2\frac{\partial^2 G}{\partial x^2}(a,X_a)\right|$$
$$\leq \mathbb{E}\left|\left[\sigma(X_u)\frac{\partial G}{\partial x}(u,X_u)\right]\left[\sigma(X_a)^2\frac{\partial^2 G}{\partial x^2}(a,X_a)\right]\left[\frac{\sigma(X_u)}{\sigma(X_a)}\right]^2\right|$$

$$+ \mathbb{E}\left|\left[\sigma(X_a)\frac{\partial G}{\partial x}(a, X_a)\right]\left[\sigma(X_a)^2\frac{\partial^2 G}{\partial x^2}(a, X_a)\right]\left[\frac{\sigma(X_u)}{\sigma(X_a)}\right]^3\right|$$

$$+ \mathbb{E}\sigma(X_a)^4 \left[\frac{\partial^2 G}{\partial x^2}(a, X_a)\right]^2 \frac{\sigma(X_u)^3 |X_u - X_a|}{\sigma(X_a)^4}.$$

Since $((\sigma\frac{\partial G}{\partial x})(u, X_u))_{u \in [0,1)}$ is an $L_2$-martingale (for a similar reason the process $((\sigma^2\frac{\partial^2 G}{\partial x^2})(u, X_u))_{u \in [0,1)}$ shares this property) we finally get that

$$\mathbb{E}\left|\Phi_i^n(u, X_u)\sigma(X_u)^2 \frac{\partial^2 G}{\partial x^2}(t_i^{n,\beta}, X_{t_i^{n,\beta}})\right| \\ \leq c_{(15)}\left[H(u)\left\|\left(\sigma\frac{\partial G}{\partial x}\right)(u, X_u)\right\|_{L_2} + H(u)^2\right] \quad (15)$$

for some absolute constant $c_{(15)} > 0$. Using $\mathbb{E}|(\sigma\frac{\partial G}{\partial x})(u, X_u)|^2 = \sum_{k=1}^{\infty} k\alpha_k^2 u^{k-1}$ for $g = \sum_{k=0}^{\infty} \alpha_k h_k$ if $X = B$ and $g(e^{\cdot - (1/2)}) = \sum_{k=0}^{\infty} \alpha_k h_k$ if $X = S$, where $(h_k)_{k=0}^{\infty}$ are the normalized Hermite polynomials and [12], Lemma 3.9, we get that

$$\left\|\left(\sigma\frac{\partial G}{\partial x}\right)(u, X_u)\right\|_{L_2} \leq c_{(16)}[\|g(X_1)\|_{L_2} + H(u)], \quad (16)$$

where $c_{(16)} > 0$ is an absolute constant, so that

$$\psi_n(u) \leq c_{(15)}\left[H(u)\left\|\left(\sigma\frac{\partial G}{\partial x}\right)(u, X_u)\right\|_{L_2} + H(u)^2\right] + c_{(14)}H(u)^2 \\ \leq [c_{(14)} + c_{(15)}]H(u)^2 + c_{(15)}c_{(16)}H(u)[\|g(X_1)\|_{L_2} + H(u)]$$

and inequality (13) follows.

(d) Now we can conclude the proof. Because of (5) the assumption of Proposition 5.1 implies that

$$\int_0^1 (1-u)^{1-\beta}|H(u) \vee \|g(X_1)\|_{L_2}|^2 \, du < \infty$$

and it remains to show that

$$\lim_n \psi_n(u) = 0 \quad \text{for all } u \in [0, 1).$$

But this follows from

$$\lim_n \mathbb{E}\sum_{i=0}^{n-1}\left|\Phi_i^n(u, X_u)\frac{\partial^2 G}{\partial x^2}(t_i^{n,\beta}, X_{t_i^{n,\beta}})\sigma(X_u)^2\right|\chi_{[t_i^{n,\beta}, t_{i+1}^{n,\beta})}(u) = 0$$

and

$$\lim_n \mathbb{E} \sum_{i=0}^{n-1} \left[ \frac{\partial^2 G}{\partial x^2}(u, X_u) - \frac{\partial^2 G}{\partial x^2}(t_i^{n,\beta}, X_{t_i^{n,\beta}}) \right]^2 \sigma(X_u)^4 \chi_{[t_i^{n,\beta}, t_{i+1}^{n,\beta}]}(u) = 0$$

by dominated convergence and Lemma 4.1. □

**Lemma 5.2.** *For $\beta \in (0,1]$ and $g(X_1) \in L_2$ one has that*

$$\lim_n n \mathbb{E}|I_1^2(\tau^{n,\beta})|^2 = \lim_n n \mathbb{E} \left| \sum_{i=0}^{n-1} \int_{t_i^{n,\beta}}^{t_{i+1}^{n,\beta}} [\sigma(X_u) - \sigma(X_{t_i^{n,\beta}})] \right.$$
$$\left. \times \frac{\partial^2 G}{\partial x^2}(t_i^{n,\beta}, X_{t_i^{n,\beta}})(X_u - X_{t_i^{n,\beta}}) \, \mathrm{d}B_u \right|^2 = 0.$$

**Proof.** For $X = B$ the integrand vanishes so that we would only need to check the case $X = S$. This can be done by observing

$$\mathbb{E} \left| \sum_{i=0}^{n-1} \int_{t_i^{n,\beta}}^{t_{i+1}^{n,\beta}} [S_u - S_{t_i^{n,\beta}}]^2 \frac{\partial^2 G}{\partial x^2}(t_i^{n,\beta}, S_{t_i^{n,\beta}}) \, \mathrm{d}B_u \right|^2$$
$$= \sum_{i=0}^{n-1} \int_{t_i^{n,\beta}}^{t_{i+1}^{n,\beta}} \mathbb{E} \left| S_{t_i^{n,\beta}}^2 \frac{\partial^2 G}{\partial x^2}(t_i^{n,\beta}, S_{t_i^{n,\beta}}) \right|^2 \mathbb{E} \left| \frac{S_u - S_{t_i^{n,\beta}}}{S_{t_i^{n,\beta}}} \right|^4 \mathrm{d}u,$$

exploiting that $H$ is continuous and non-decreasing, and by using

$$\int_0^1 (1-u) H(u)^2 \, \mathrm{d}u < \infty,$$

which is true for all $g$ with $g(X_1) \in L_2$ and a consequence of [12], Lemma 3.9 (cf. (5)). □

*Preparations for $I^3(\tau)$.*

The process $I^3(\tau)$ is responsible for the structure of the weak limit of the renormalized error process. The next lemma is a counterpart of [22], Lemma 1.5. For the convenience of the reader we give some details concerning the proof in the Appendix.

**Lemma 5.3.** *Let $k \in \{1, 2\}$, $T \in (0, 1]$ and let $a = (a_t)_{t \in [0,T]}$ be a continuous process. Define*

$$\psi_s^{n,k}(a) := n^{k/2} \sum_{i=0}^{n-1} a_{t_i^{n,\beta}} \left( \frac{X_s - X_{t_i^{n,\beta}}}{\sigma(X_{t_i^{n,\beta}})} \right)^k \chi_{[t_i^{n,\beta}, t_{i+1}^{n,\beta})}(s)$$

for $s \in [0, T]$. Then

$$\lim_n \left[ \sup_{t \in [0,T]} \left| \int_0^t \psi_s^{n,1}(a) \, ds \right| + \sup_{t \in [0,T]} \left| \int_0^t \psi_s^{n,2}(a) \, ds - \frac{1}{2} \int_0^t \nu_\beta(s) a_s \, ds \right| \right] = 0$$

in probability, where $\nu_\beta(s) = (1/\beta)(1-s)^{1-\beta}$.

**Lemma 5.4.** *For $T \in (0,1)$ one has that*

$$\left\| \sup_{t \in [T,1]} |\sqrt{n} C_t(\tau^{n,\beta}) - \sqrt{n} C_T(\tau^{n,\beta})| \right\|_{L_2}$$

$$\leq \frac{c}{\sqrt{\beta}} \left( \int_{(T-1/(\beta n))^+}^1 (1-s)^{1-\beta} H(s)^2 \, ds \right)^{1/2},$$

*where $H(s)^2 = \mathbb{E}((\sigma^2 \frac{\partial^2 G}{\partial x^2})(s, X_s))^2$ and $c > 0$ is an absolute constant.*

**Proof.** Let $T_n^\beta := \sup\{t_i^{n,\beta} : t_i^{n,\beta} \leq T, i = 0, \ldots, n-1\}$. Then, by Doob's maximal inequality and [10], Proof of Theorem 6,

$$\left\| \sup_{t \in [T,1]} |C_t(\tau^{n,\beta}) - C_T(\tau^{n,\beta})| \right\|_{L_2} \leq 4 \| C_1(\tau^{n,\beta}) - C_{T_n^\beta}(\tau^{n,\beta}) \|_{L_2}$$

$$\leq c \left( \sum_{i=0}^{n-1} \int_{t_i^{n,\beta} \vee T_n^\beta}^{t_{i+1}^{n,\beta} \vee T_n^\beta} (t_{i+1}^{n,\beta} - s) H(s)^2 \, ds \right)^{1/2}$$

$$\leq c \sup_{\substack{i=0,\ldots,n-1 \\ s \in [t_i^{n,\beta}, t_{i+1}^{n,\beta})}} \left| \frac{t_{i+1}^{n,\beta} - s}{(1-s)^{1-\beta}} \right|^{1/2} \left( \int_{T_n^\beta}^1 (1-s)^{1-\beta} H(s)^2 \, ds \right)^{1/2}$$

$$\leq \frac{c}{\sqrt{\beta n}} \left( \int_{(T-1/(\beta n))^+}^1 (1-s)^{1-\beta} H(s)^2 \, ds \right)^{1/2},$$

where $c > 0$ is an absolute constant and we have used (4). □

The next theorem is due to Rootzén and was formulated for $T = 1$.

**Theorem 5.5 ([[22]], Theorem 1.2).** *Let $T \in [0,1]$. Suppose that $\psi^n = (\psi_t^n)_{t \in [0,T]}$, $n = 1, 2, \ldots$, are progressively measurable processes and that*

$$\lim_n \sup_{t \in [0,T]} \left| \int_0^t [\psi_s^n]^2 \, ds - A_t \right| = 0 \qquad \text{in probability}$$

*for some continuous process $A = (A_t)_{t \in [0,T]}$ and that*

$$\lim_n \sup_{t \in [0,T]} \left| \int_0^t \psi_s^n \, ds \right| = 0 \qquad \text{in probability.}$$

*Then*

$$\left( \int_0^t \psi_s^n \, dB_s \right)_{t \in [0,T]} \Longrightarrow_{C[0,T]} (W_{A_t})_{t \in [0,T]} \qquad \text{for } n \to \infty,$$

*where the Brownian motion $W$ is independent from $\mathcal{F}$.*

**Proof of Theorem 3.1.** Combining Theorem 5.5 and Lemma 5.3 for $a_t := (\sigma^2 \frac{\partial^2 G}{\partial x^2})(t, X_t)$ in case $k = 1$, $a_t := [(\sigma^2 \frac{\partial^2 G}{\partial x^2})(t, X_t)]^2$ in case $k = 2$ and $A_t := A_\beta(t)$ yields to

$$(\sqrt{n} I_t^3(\tau^{n,\beta}))_{t \in [0,T]} \Longrightarrow_{C[0,T]} (W_{A_\beta(t)})_{t \in [0,T]} = (Z_\beta(t))_{t \in [0,T]}$$

for all $T \in [0,1)$. Because of Proposition 5.1, Lemma 5.2 and Doob's maximal inequality (note that $(I_t^1(\tau^{n,\beta}))_{t \in [0,T]}$ and $(I_t^2(\tau^{n,\beta}))_{t \in [0,T]}$ are $L_2$-martingales (cf. Lemma 4.1) so that $\sqrt{n} \sup_{t \in [0,T]} |I_t^k(\tau^{n,\beta})| \to_{L_2} 0$ as $n \to \infty$ for $k = 1, 2$), we can deduce that

$$(\sqrt{n} C_t(\tau^{n,\beta}))_{t \in [0,T]} \Longrightarrow_{C[0,T]} (Z_\beta(t))_{t \in [0,T]} \qquad \text{as } n \to \infty. \tag{17}$$

*Proof of* (i) $\Rightarrow$ (iii): First, we observe that (i) implies

$$\int_0^1 (1-s)^{1-\beta} H(s)^2 \, ds < \infty$$

according to (5) so that $\mathbb{E} A_\beta(1) < \infty$. Given a continuous and bounded $\varphi : C[0,1] \to \mathbb{R}$ we have to prove that

$$\lim_n \mathbb{E} \varphi(Y^n) = \mathbb{E} \varphi(\widetilde{Z}_\beta),$$

where $Y_t^n := \sqrt{n} C_t(\tau^{n,\beta})$. We can restrict ourselves to uniformly continuous and bounded $\varphi$ (cf. [3]). Let $T \in (0,1)$, $Y^{T,n} := (Y_{t \wedge T}^n)_{t \in [0,1]}$, and $\widetilde{Z}_\beta^T := (\widetilde{Z}_\beta(t \wedge T))_{t \in [0,1]}$. Then

$$|\mathbb{E}\varphi(Y^n) - \mathbb{E}\varphi(\widetilde{Z}_\beta)| \leq |\mathbb{E}\varphi(Y^n) - \mathbb{E}\varphi(Y^{T,n})| + |\mathbb{E}\varphi(Y^{T,n}) - \mathbb{E}\varphi(\widetilde{Z}_\beta^T)|$$
$$+ |\mathbb{E}\varphi(\widetilde{Z}_\beta^T) - \mathbb{E}\varphi(\widetilde{Z}_\beta)|.$$

We fix $\varepsilon > 0$ and find a $\delta > 0$ such that $|\varphi(f) - \varphi(g)| < \varepsilon$ for $\|f - g\|_{C[0,1]} < \delta$. Then

$$|\mathbb{E}\varphi(Y^n) - \mathbb{E}\varphi(Y^{T,n})|$$
$$\leq \int_{\{\|Y^n - Y^{T,n}\|_{C[0,1]} \geq \delta\}} |\varphi(Y^n) - \varphi(Y^{T,n})| \, d\mathbb{P}$$

$$+ \int_{\{\|Y^n - Y^{T,n}\|_{C[0,1]} < \delta\}} |\varphi(Y^n) - \varphi(Y^{T,n})| \, d\mathbb{P}$$

$$\leq 2\|\varphi\|_\infty \mathbb{P}(\|Y^n - Y^{T,n}\|_{C[0,1]} \geq \delta) + \varepsilon$$

$$\leq 2\|\varphi\|_\infty \frac{c_{(5.4)}^2}{\delta^2 \beta} \int_{(T - 1/(\beta n))^+}^1 (1-s)^{1-\beta} H(s)^2 \, ds + \varepsilon,$$

where we have used Lemma 5.4. Let $T_0 \in (0,1)$ be such that

$$2\|\varphi\|_\infty \frac{c_{(5.4)}^2}{\delta^2 \beta} \int_{T_0}^1 (1-s)^{1-\beta} H(s)^2 \, ds \leq \varepsilon$$

and

$$|\mathbb{E}\varphi(\widetilde{Z}_\beta^T) - \mathbb{E}\varphi(\widetilde{Z}_\beta)| \leq \varepsilon$$

for $T \in [T_0, 1)$ (note that $\|\widetilde{Z}_\beta^T(\omega) - \widetilde{Z}_\beta(\omega)\|_{C[0,1]} \to 0$ as $T \uparrow 1$ for all $\omega \in \Omega$). Fix $n_0 \geq 1$ such that $1/(\beta n_0) \leq (1 - T_0)/2$. Hence, for $T \in [(T_0 + 1)/2, 1)$ and $n \geq n_0$ one has $T - \frac{1}{\beta n} \geq T_0$ and

$$|\mathbb{E}\varphi(Y^n) - \mathbb{E}\varphi(\widetilde{Z}_\beta)| \leq 3\varepsilon + |\mathbb{E}\varphi(Y^{T,n}) - \mathbb{E}\varphi(\widetilde{Z}_\beta^T)|.$$

Defining the bounded and continuous function $\varphi_T : C[0,T] \to \mathbb{R}$ by $\varphi_T(g) := \varphi(f)$ with $f(t) := g(t \wedge T)$, we get

$$\lim_n \mathbb{E}\varphi(Y^{T,n}) = \lim_n \mathbb{E}\varphi_T((Y_t^n)_{t \in [0,T]}) = \mathbb{E}\varphi_T((\widetilde{Z}_\beta(t))_{t \in [0,T]}) = \mathbb{E}\varphi(\widetilde{Z}_\beta^T),$$

where we used (17), and

$$\limsup_n |\mathbb{E}\varphi(Y^n) - \mathbb{E}\varphi(\widetilde{Z}_\beta)| \leq 3\varepsilon.$$

Since this is true for all $\varepsilon > 0$ we are done. $\square$

# Appendix

First we formalize some ideas of [22].

**Lemma A.1.** *Let $T \in (0,1]$ and $\mu_n(\omega) = \mu_n^+(\omega) - \mu_n^-(\omega)$, where $\mu_n^+(\omega)$ and $\mu_n^-(\omega)$ are finite Borel measures on $[0,T]$ for $\omega \in \Omega$. Assume that*

  (i) *$\mu_n^\pm([0,t])$ are measurable for all $t \in [0,T]$,*
  (ii) *$\sup_n \mathbb{E}|(\mu_n^+ + \mu_n^-)([0,T])|^p < \infty$ for some $p \in (0, \infty)$ and*

(iii) *there is a finite Borel measure $\mu$ on $[0,T]$ such that, in probability,*

$$\lim_n \sup_{t\in[0,T]} |\mu_n([0,t]) - \mu([0,t])| = 0.$$

*Then, given a continuous process $(a_s)_{s\in[0,T]}$ of $\mathcal{F}$-measurable random variables, one has that*

$$\lim_n \sup_{t\in[0,T]} \left| \int_{[0,t]} a_s \, \mathrm{d}\mu_n(s) - \int_{[0,t]} a_s \, \mathrm{d}\mu(s) \right| = 0 \qquad \text{in probability.} \tag{18}$$

**Proof.** (a) For $N = 1, 2, \ldots$ let

$$a_t^N := a_0 \chi_{[0,T/2^N]}(t) + \sum_{l=1}^{2^N-1} a_{(l/2^N)T} \chi_{((l/2^N)T, ((l+1)/2^N)T]}(t)$$

$$= a_{T(2^N-1)/2^N} \chi_{[0,T]}(t) + (a_{T(2^N-2)/2^N} - a_{T(2^N-1)/2^N}) \chi_{[0,T(2^N-1)/2^N]}(t) + \cdots$$

$$+ (a_0 - a_{T/2^N}) \chi_{[0,T1/2^N]}(t).$$

To check (18) for $a^N = (a_t^N)_{t\in[0,T]}$ it is enough to verify (18) for $a^N$ replaced by $b = (\varphi \chi_{[0,r]}(t))_{t\in[0,T]}$ with $r \in [0,T]$ and an $\mathcal{F}$-measurable random variable $\varphi$. Since $\varphi$ is a constant factor, an easy argument shows that it is sufficient to check the case $\varphi \equiv 1$. But then we can use Assumption (iii) and obtain (18) for $a^N$.

(b) To replace $a^N$ by $a$ we observe that

$$\sup_{t\in[0,T]} \left| \int_{[0,t]} a_s \, \mathrm{d}\mu_n(s) - \int_{[0,t]} a_s \, \mathrm{d}\mu(s) \right| \leq \sup_{t\in[0,T]} |a_t - a_t^N|(\mu_n^+ + \mu_n^- + \mu)([0,T])$$

$$+ \sup_{t\in[0,T]} \left| \int_{[0,t]} a_s^N \, \mathrm{d}\mu_n(s) - \int_{[0,t]} a_s^N \, \mathrm{d}\mu(s) \right|.$$

Because of (ii) and $\sup_{t\in[0,T]} |a_t(\omega) - a_t^N(\omega)| \to 0$ as $N \to \infty$ for all $\omega \in \Omega$ step (a) implies the assertion. $\square$

**Proof of Lemma 5.3.** The proof is similar to the one in [22]; the part that differs is the estimate of (22). For the convenience of the reader we give some details. Define the random measures

$$\mu_n^k := n^{k/2} \sum_{i=0}^{n-1} \delta_{\{s_i^{n,\beta}\}} \int_{s_i^{n,\beta}}^{s_{i+1}^{n,\beta}} \left( \frac{X_s - X_{s_i^{n,\beta}}}{\sigma(X_{s_i^{n,\beta}})} \right)^k \mathrm{d}s$$

for $k \in \{1,2\}$ and $s_i^{n,\beta} := t_i^{n,\beta} \wedge T$ and let $(\mu_n^k)^\pm$ be the positive and negative parts ($\omega$-wise), respectively. By a standard computation one checks that

$$\sup_n \mathbb{E} n^{k/2} \sum_{i=0}^{n-1} \int_{s_i^{n,\beta}}^{s_{i+1}^{n,\beta}} \left| \frac{X_s - X_{s_i^{n,\beta}}}{\sigma(X_{s_i^{n,\beta}})} \right|^k \mathrm{d}s < \infty$$

so that $\sup_n \mathbb{E}((\mu_n^k)^+ + (\mu_n^k)^-)([0,T]) < \infty$. Moreover, using (4),

$$\sup_{t\in[0,T]} \left| \int_0^t \psi_s^{n,k}(a)\,\mathrm{d}s - \int_{[0,t]} a_s\,\mathrm{d}\mu_n^k(s) \right|^2 \leq \sup_{\substack{0\leq i\leq n-1 \\ t\in[s_i^{n,\beta}, s_{i+1}^{n,\beta}]}} \left| \int_t^{s_{i+1}^{n,\beta}} \psi_s^{n,k}(a)\,\mathrm{d}s \right|^2$$

$$\leq \sup_{0\leq i\leq n-1} (s_{i+1}^{n,\beta} - s_i^{n,\beta}) \int_0^T (\psi_s^{n,k}(a))^2\,\mathrm{d}s$$

$$\leq \frac{(a^*)^2}{\beta n} \int_0^T (\psi_s^{n,k}(1))^2\,\mathrm{d}s,$$

where $a^* := \sup_{t\in[0,T]} |a_t|$ and $\sup_{m\geq 1} \mathbb{E}\int_0^T (\psi_s^{m,k}(1))^2\,\mathrm{d}s < \infty$, so that

$$\lim_n \sup_{t\in[0,T]} \left| \int_0^t \psi_s^{n,k}(a)\,\mathrm{d}s - \int_{[0,t]} a_s\,\mathrm{d}\mu_n^k(s) \right| = 0 \tag{19}$$

in probability. In view of (19) and Lemma A.1 we only need to verify

$$\lim_n \sup_{t\in[0,T]} |\mu_n^1([0,t])| = 0 \quad \text{and} \quad \lim_n \sup_{t\in[0,T]} \left| \mu_n^2([0,t]) - \frac{1}{2}\int_0^t \nu_\beta(s)\,\mathrm{d}s \right| = 0 \tag{20}$$

in probability. Let $E_1 := 0$, $E_2 := 1/2$ and $b_s^{n,k} := \psi_s^{n,k}(1)$. To show (20) we upper bound

$$\sup_{0\leq i\leq n} \left| \int_0^{s_i^{n,\beta}} b_s^{n,k}\,\mathrm{d}s - E_k \int_0^{s_i^{n,\beta}} \nu_\beta(s)\,\mathrm{d}s \right| + \sup_{1\leq i\leq n} \sup_{t\in[s_{i-1}^{n,\beta}, s_i^{n,\beta}]} \left| E_k \int_t^{s_i^{n,\beta}} \nu_\beta(s)\,\mathrm{d}s \right|$$

$$\leq \sup_{0\leq i\leq n} \left| \int_0^{s_i^{n,\beta}} b_s^{n,k}\,\mathrm{d}s - \mathbb{E}\int_0^{s_i^{n,\beta}} b_s^{n,k}\,\mathrm{d}s \right| \tag{21}$$

$$+ \sup_{0\leq i\leq n} \left| \mathbb{E}\int_0^{s_i^{n,\beta}} b_s^{n,k}\,\mathrm{d}s - E_k \int_0^{s_i^{n,\beta}} \nu_\beta(s)\,\mathrm{d}s \right| + 1/(2\beta^2 n), \tag{22}$$

where we used (4) again.

Term (21): By Doob's maximal inequality for martingales one can show that

$$\mathbb{E}\sup_{0\leq i\leq n} \left| \int_0^{s_i^{n,\beta}} b_s^{n,k}\,\mathrm{d}s - \mathbb{E}\int_0^{s_i^{n,\beta}} b_s^{n,k}\,\mathrm{d}s \right|^2 \xrightarrow[n\to\infty]{} 0.$$

Term (22): Because for $k=1$ the term is zero we assume that $k=2$ and get

$$\sup_{0\leq i\leq n} \left| \mathbb{E}\int_0^{s_i^{n,\beta}} b_s^{n,2}\,\mathrm{d}s - \frac{1}{2}\int_0^{s_i^{n,\beta}} \nu_\beta(s)\,\mathrm{d}s \right|$$

$$= \sup_{1\leq i\leq n} \left| \sum_{j=0}^{i-1} \int_{s_j^{n,\beta}}^{s_{j+1}^{n,\beta}} n\mathbb{E}\left(\frac{X_s - X_{s_j^{n,\beta}}}{\sigma(X_{s_j^{n,\beta}})}\right)^2 \mathrm{d}s - \frac{1}{2}\int_0^{s_i^{n,\beta}} \nu_\beta(s)\,\mathrm{d}s \right|$$

$$= \sup_{1\leq i\leq n} \left| \sum_{j=0}^{i-1} \int_{s_j^{n,\beta}}^{s_{j+1}^{n,\beta}} n(s - s_j^{n,\beta} + m(s - s_j^{n,\beta})(s - s_j^{n,\beta})^2)\,\mathrm{d}s - \frac{1}{2}\int_0^{s_i^{n,\beta}} \nu_\beta(s)\,\mathrm{d}s \right|$$

$$\leq \sum_{j=0}^{n-1} \int_{s_j^{n,\beta}}^{s_{j+1}^{n,\beta}} nm(s - s_j^{n,\beta})(s - s_j^{n,\beta})^2\,\mathrm{d}s$$

$$+ \frac{1}{2}\sum_{j=0}^{n-1} \int_{s_j^{n,\beta}}^{s_{j+1}^{n,\beta}} |n(s_{j+1}^{n,\beta} - s_j^{n,\beta}) - \nu_\beta(s)|\,\mathrm{d}s$$

$$\leq \frac{\mathrm{e}}{3\beta^3 n} + \frac{1}{2}\sum_{j=0}^{n-1} \int_{s_j^{n,\beta}}^{s_{j+1}^{n,\beta}} |n(s_{j+1}^{n,\beta} - s_j^{n,\beta}) - \nu_\beta(s)|\,\mathrm{d}s$$

$$\leq \frac{\mathrm{e}}{3\beta^3 n} + \frac{1}{2}\int_{[0,t_{j_{0,n}}^{\beta,n})} \frac{1}{\beta} \sup_{u\in[0,1-1/n]} \left|\left(u + \frac{1}{n}\right)^{1/\beta-1} - u^{1/\beta-1}\right|\,\mathrm{d}s$$

$$+ \frac{1}{2}\int_{[t_{j_{0,n}}^{\beta,n},T]} |n(s_{j+1}^{n,\beta} - s_j^{n,\beta}) - \nu_\beta(s)|\,\mathrm{d}s,$$

where $m\colon [0,1] \to [0,\mathrm{e}]$ is a continuous function and $j_{0,n}$ is the largest $j \in \{0,1,\ldots,n\}$ such that $s_j^{n,\beta} = t_j^{n,\beta}$. Finally, by (4) we can bound the last term by

$$\frac{1}{\beta}|T - t_{j_{0,n}}^{n,\beta}| \leq \frac{1}{\beta^2 n}$$

so that the term (22) converges to zero as $n \to \infty$ and the proof is complete. $\square$

## Acknowledgement

We would like to thank the referee for the careful reading of the manuscript and his hints that improved the readability of the paper.